# Continuity in Discrete Sets

**Mark Burgin**

Department of Mathematics, UCLA
405 Hilgard Ave.
Los Angeles, 90046, USA

**Abstract**: Continuous models used in physics and other areas of mathematics applications become discrete when they are computerized, e.g., utilized for computations. Besides, computers are controlling processes in discrete spaces, such as films and television programs. At the same time, continuous models that are in the background of discrete representations use mathematical technology developed for continuous media. The most important example of such a technology is calculus, which is so useful in physics and other sciences. The main goal of this paper is to synthesize continuous features and powerful technology of the classical calculus with the discrete approach of numerical mathematics and computational physics. To do this, we further develop the theory of fuzzy continuous functions and apply this theory to functions defined on discrete sets. The main interest is the classical Intermediate Value theorem. Although the result of this theorem is completely based on continuity, utilization of a relaxed version of continuity called fuzzy continuity, allows us to prove discrete versions of the Intermediate Value theorem. This result provides foundations for a new approach to discrete dynamics.

Contemporary science and engineering are based on mathematical tools and methods. Calculus is one of the most efficient tools for modeling real world phenomena. That is why calculus is so useful in science and engineering. However, in its essence calculus is opposite to real world because calculus is exact, rigorous and abstract while real things and systems are imprecise, vague, and concrete. To lessen this gap, mathematicians create new methods, making possible to work with natural vagueness and incompleteness of information using exact mathematical structures.

The necessity to launch investigation and implementation of fuzzy principles in the classical analysis and calculus, while studying ordinary functions, is caused by several reasons. One of the most important of them is connected to properties of measurements. Any real measurement does not provide absolutely precise results but only approximate values. For example, it is impossible to find out if any series of numbers obtained in experiments converges or a function determined by measurements is continuous at a given point. Consequently, constructions and methods developed in the classical analysis are only approximations to what exists in reality. In many situations, such approximations has been giving a sufficiently adequate representation of studied phenomena. However, scientists and, especially, engineers have discovered many cases in which such methods did not work because the classical approach is too rough.

For instance, there are three common problems in modeling (1) identifying the "correct" model, (2) estimating parameters from experimental data, and (3) determining which of many parameters contribute most to the uncertainties in the prediction of the model. The issue is not if there are uncertainties, they are always present, the real challenge is to identify those components of the model that have the most influence on the prediction.

This approach correlates with problems of modern physical theories in which physical systems are described by chaotic processes. Taking into account the fact that chaotic solutions are obtained by computations, physicists ask (Cartwrite and Piro, 1992; Gontar, 1997) whether chaotic solutions of the differential equations that model physical systems reflect the dynamic laws of nature represented by these equations or whether they are solely the result of an extreme sensitivity of these solutions to numerical procedures and computational errors. To solve these and other problems of the theory of chaos, we need adequate and flexible calculus on discrete sets.

Consequently, new methods and constructions are necessary to take into account such more sophisticated effects that exist in different systems. Such methods and constructions are provided by neoclassical analysis (cf., (Burgin, 1995; 1999; 2008)). Neoclassical analysis extends the scope of analysis making, at the same time, its methods more precise. Consequently, new results are obtained extending and even completing classical theorems. In addition, facilities of analytical methods for various applications also become broader and more efficient.

In this paper, we apply neoclassical analysis to study fuzzy continuous functions. The Intermediate Value Theorem is one of the central classical results for continuous functions. It is important for many problems of calculus and its applications.

Let us consider a continuous function $f: [a, b] \to \mathbf{R}$.

**Intermediate Value Theorem.** If $f(a) = k$, $f(b) = h$, and $k < h$ ( $k > h$), then for any number $l$ from the interval $[k, h]$ ( from the interval $[h, k]$ ) there is [at least one] point $c \in [a, b]$ such that $f(c) = l$.

When the function $f(x)$ is not continuous or it is defined not for all points in the interval $[a, b]$, then this result is not true in a general case. Simple examples show that when a function has gaps, for example, when it is defined on a discrete set, then this result fails.

However, fuzzy continuity studied in (Burgin, 1995; 1999; 2008) allows us to find similar properties for functions defined not for all points in the interval or that are not continuous in the interval.

Main constructions and concepts from neoclassical analysis are introduced in Section 2. In Section 3, several extensions of the classical Intermediate Value Theorem (the Fuzzy Intermediate Value Theorem, Discrete Intermediate Value Theorem and Digital Topology Intermediate Value Theorem) are proved. In Section 4, these theorems are applied to study images of connected and fuzzy connected sets, as well as to prove existence and fuzzy continuity of inverse functions. Results for discrete continuous functions are applied to functions in digital topology (cf. (Rosenfeld, 1979; 1986; Boxer, 1994)).

It is necessary to remark that application of neoclassical analysis to functions on discrete sets extends methods of digital topology, making them efficient in situations that exist in image recognition and processing where digital topology cannot give sufficiently precise results. The problem is that digital topology works with uniform grids and all its principal constructions, such as continuity and connectedness, are based on utilization of the minimal distance between points in the grid, e.g., pixels in the screen of a TV set or computer monitor. However, many computational problems demand utilization of a non-uniform mesh or grid. For instance, to numerically solve complex set of PDEs with solutions varying on different length scales, it is necessary to use advanced computational methods and perform the calculations on high performance computers, typically in parallel on many processors. The single most important factor is to use adaptive grid refinement. In some cases, it useful when the grid has the hierarchical structure that would adapt according to some preset criteria (cf., for example, (Hänel, *et al*, 1999)).

Besides, even when the grid is uniform, as the grid of pixels is, it might be necessary to treat different parts of the grid with a specific precision. For instance, at first, the coarsest grid (also known as the base grid) covers the whole computational domain. In places where more resolution is needed a finer (child) grid is created on top of the coarser one (the parent grid). These grids are layered on top of each other in hierarchical fashion, i.e., each child grid is fully contained within a parent grid. The definition of grid refinement is, in some cases, motivated by the requirement that grid point concentration should be higher at locations of rapid variation in the function values.

This implies necessity in scalable calculus constructions. Such constructions are developed in this paper. Scalability is achieved by variation of parameters $r$, $q$ and others because these parameters determine precision of continuity, connectedness, differentiation, integration and other topological and analytical constructions and tools.

**2. Fuzzy Continuous Functions**

We start with describing some structures introduced earlier and introducing some new concepts to the neoclassical analysis.

Let $\mathbf{R}$ be a set of all real numbers, $\mathbf{R}^+$ be a set of all non-negative real numbers, and $r \in \mathbf{R}^+$.

**Definition 1.** A function $f: \mathbf{R} \to \mathbf{R}$ is called *r-continuous* at a point $a \in \mathbf{R}$ if $f(x)$ is defined at $a$ and for any $\varepsilon > 0$ there is $\delta > 0$ such that the inequality $|a - x| < \delta$ implies the inequality $|f(x) - f(a)| < r + \varepsilon$, or in other words, for any $x$ with $|a - x| < \delta$, we have $|f(x) - f(a)| < r + \varepsilon$.

**Example 1.** Let us consider the function $f(x)$ equal to $x^3$ when $x \in [0, 1/2]$ and equal to $x$ at all other point of the real line. This function is (3/8)-continuous at the point ½ and 0-continuous at all other points.

**Definition 2.** A function $f: \mathbf{R} \to \mathbf{R}$ is called *r-continuous in* (*inside*) $X$ if $f(x)$ (the restriction of $f(x)$ on $X$) is *r*-continuous at each point $a$ from $X \cap \text{Dom } f$.

**Example 2.** The function $f(x)$ from Example 1 is (3/8)-continuous on the whole real line $\mathbf{R}$.

**Lemma 1** (Burgin, 2008). A function $f(x)$ is 0-continuous at a point $a \in \mathbf{R}$ if and only if it is continuous at the point $a$.

**Corollary 1.** A function $f(x)$ is 0-continuous if and only if it is continuous.

These results show that the concept of $r$-continuity is a natural extension of the concept of conventional continuity.

**Definition 3.** A function $f: \mathbf{R} \to \mathbf{R}$ is called $(q, r)$-*continuous at a point* $a \in \mathbf{R}$ if for any $\varepsilon > 0$ there is $\delta > 0$ such that the inequality $|a - x| < q + \delta$ implies the inequality $|f(x) - f(a)| < r + \varepsilon$, or in other words, for any $x$ with $|a - x| < q + \delta$, we have $|f(x) - f(a)| < r + \varepsilon$.

**Example 3.** Let us consider functions $f(x) = x^3$ and $g(x) = x$ in the interval $[0, 2]$. The function $g(x)$ is $(0.1, 0.1)$-continuous at the point 1, while the function $f(x)$ is not $(0.1, 0.1)$-continuous at this point. The function $f(x)$ is $(0.03, 0.1)$-continuous at this point.

**Definition 4.** A function $f: \mathbf{R} \to \mathbf{R}$ is called $(q, r)$-*continuous in* (*inside*) $X$ if $f(x)$ (the restriction of $f(x)$ on $X$) is $(q, r)$-continuous at each point $a$ from $X \cap \mathrm{Dom}\, f$.

**Lemma 2** (Burgin, 2008). a) A function $f(x)$ is $r$-continuous (at a point $a \in \mathbf{R}$) if and only if it is $(0, r)$-continuous (at the point $a$).

b) A function $f(x)$ is continuous (at a point $a \in \mathbf{R}$) if and only if it is $(0,0)$-continuous (at the point $a$).

These results show that the concept of $(q, r)$-continuity is a natural extension of the concept of conventional continuity.

Definitions 3 and 4 imply the following result.

**Lemma 3** (Burgin, 2008). If $t > r$, and $p < q$, then any $(q, r)$-continuous (at $a$) function $f(x)$ is also $(p, t)$-continuous (at $a$).

Indeed, if that the inequality $|a - x| < q + \delta$ implies the inequality $|f(x) - f(a)| < r + \varepsilon$, then the inequality $|a - x| < p + \delta$ implies the inequality $|a - x| < q + \delta$ because $p < q$. The inequality $|a - x| < q + \delta$ implies the inequality $|f(x) - f(a)| < r + \varepsilon$ as $f(x)$ is $(q, r)$-continuous at $a$, which, in turn, implies the inequality $|f(x) - f(a)| < t + \varepsilon$ because $t > r$. Thus, by Definition 3, if $f(x)$ is $(q, r)$-continuous at a point $a$ function, then it is also $(p, t)$-continuous at $a$.

**Corollary 2.** If $q > l$ ($r < p$), then any $(q, r)$-continuous at a point $a \in R$ function is also $(l, r)$-continuous (also $(q, p)$-continuous) at the point $a$.

**Corollary 3.** If $q > l$ and $r < p$, then any $(q, r)$-continuous function is also $(l, p)$-continuous.

**Corollary 4.** If a function $f(x)$ is $(q, r)$-continuous at a point $a \in R$, then $f(x)$ is $r$-continuous at the point $a$.

Let us consider a subset $X$ of the real line $R$.

**Definition 5.** a) The set $X$ is *discrete* in $R$ if there is no sequence $l = \{ a_i \in X; i = 1, 2, 3, \ldots \}$ such that $a = \lim l$ for some real number $a$.

b) The set $X$ is *discrete* in itself if there is no sequence $l = \{ a_i \in X; i = 1, 2, 3, \ldots \}$ such that $a = \lim l$ for some real number $a \in X$.

**Example 4.** Let $u \in R$. Sets $X = \{ ku; k = 0, \pm 1, \pm 2, \pm 3, \ldots \}$ and $Z = \{ ku; k \in Z$ and $m < k < n \}$ are discrete in $R$.

**Example 5.** The set $X = \{ 0, 1, 2, 2\frac{1}{2}, 3, 3_{1/3}, 3_{2/3}, 4, 4\frac{1}{4}, 4\frac{1}{2}, 4\frac{3}{4}, 5, \ldots \}$ is discrete in $R$.

**Example 6.** The set $X = \{ 1/(2^n); n = 0, 1, 2, 3, \ldots \}$ is discrete in itself but it is not discrete in $R$.

**Remark 1.** The image of a discrete set can be not a discrete set in a general case. For instance, the set $W = \{ 0, 1, 2, 3, \ldots \}$ of all whole numbers is discrete, while the set $X = \{ 0, 1, \frac{1}{2}, 1/3, \frac{1}{4}, 1/5, \ldots \}$ is not discrete. At the same time the function $f(x)$ defined as $f(0) = 0$ and $f(n) = 1/n$ for all $n = 1, 2, 3, \ldots$ maps $W$ onto $X$.

**Lemma 4.** If a set $X$ is discrete in $R$, it is also discrete in itself.

There are other conditions that are equivalent to discreteness of real number sets.

**Proposition 1.** The following conditions are equivalent:

a) The set $X$ is discrete in $R$.

b) In any finite interval $I$, there are only a finite number of points from $X$.

c) For any finite interval $I$, there is a number $q > 0$ such that the distance between any two points from $X \cap I$ is larger than $q$.

Proof. **(a)** ⇒ **(b)**. Let us take a discrete in **R** set X and assume that there is a finite interval **I** that contains an infinite number of points from X. Because the interval **I** is a compact set, there is a sequence $l$ = { $a_i$ ∈ X; $i$ = 1, 2, 3, … } such that $a$ = lim $l$ for some real number $a$ (Kuratowski, 1966). This contradicts to discreteness of X and by the Principle of Excluded Middle, (a) ⇒ (b).

**(b)** ⇒ **(c)** because for any finite set in **R**, there is a number $q$ > 0 such that the distance between any two points from this set is larger than $q$.

**(c)** ⇒ **(a)**. Let us take a set X in **R** such that for any finite interval **I**, there is a number $q$ > 0 such that the distance between any two points from X ∩ **I** is larger than $q$ and assume that the set X is not discrete. Then by Definition 5a, there is a sequence $l$ = { $a_i$ ∈ X; $i$ = 1, 2, 3, … } such that $a$ = lim $l$ for some real number $a$. By the definition of a limit (Ross, 1996), in any finite interval **I** that contains $a$, there are infinitely many points from X. Taking smaller and smaller interval that contain $a$, we see that for any number $q$ > 0, there are always two points from X ∩ **I** such that distance between them is smaller than $q$. This contradicts to the choice of X and by the Principle of Excluded Middle, (c) ⇒ (a).

Proposition is proved because logical implication is transitive.

**Lemma 5.** Any discrete set X in **R** is countable.

Indeed, by Proposition 1, if a set X is discrete in **R**, then any finite interval **I** contains only a finite number of points from X. As it is possible to represent **R** as the union of countable finite intervals each of which contains a finite number of points from X, the set X is countable.

**Proposition 2.** Any subset Y of a discrete in **R** (in itself) set X is discrete in **R** (in itself).

**Proposition 3.** If X and Y are discrete in **R** (in itself) sets, then X ∪ Y, X ∩ Y, and X\Y, are also discrete in **R** (in itself).

**Definition 6.** Points $a$ and $b$ from a discrete set X in **R** are called *adjacent*, or *consecutive*, if there no other points from X between them.

**Definition 7.** A discrete set X in **R** is called:

1) *uniform* if all distances between any two adjacent points in X are equal;

2) *bounded down* if all distances between any two adjacent points in X are larger than some number $k > 0$;

3) *bounded up* if all distances between any two adjacent points in X are less than some number $k > 0$.

**Example 7.** Let $u \in R$. The set $X = \{ ku; k = 0, \pm 1, \pm 2, \pm 3, \ldots \}$ is an infinite *uniform* discrete set in $R$ and the set $Z = \{ ku; k \in Z \text{ and } m < k < n \}$ is a finite uniform discrete set in $R$.

**Example 8.** The set $X = \{ 0, 1, 2, 2\frac{1}{2}, 3, 3_{1/3}, 3_{2/3}, 4, 4\frac{1}{4}, 4\frac{1}{2}, 4\frac{3}{4}, 5, \ldots \}$ is a non-uniform discrete set in $R$.

**Example 9.** The discrete set $X = \{ n^2 ; n = 1, 2, 3, \ldots \}$ is unbounded up in $R$. It is also a non-uniform discrete set.

**Proposition 4.** Any subset Y of a bounded down discrete set X is discrete is bounded down.

**Remark 2.** It is possible that a subset of a uniform discrete set is not uniform and that a subset of a bounded up discrete set X is not bounded up.

Proposition 1 implies the following results.

**Corollary 5.** In a finite interval $I$, any discrete set X is bounded up and down.

**Corollary 6.** Any finite discrete set X is bounded up and down.

**Definition 8.** 1) In a uniform discrete set (grid) X in $R$, the distances between any two consecutive points is called the *spacing* of X.

2) In a bounded down discrete set X in $R$, the lower bound of distances between adjacent points in X is called the *lower inner bound* of X; it is denoted by libX.

3) In a bounded up discrete set X in $R$, the lower bound of distances between adjacent points in X is called the *upper inner bound* of X; it is denoted by uibX.

**Proposition 5.** a) libX $\leq$ uibX. b) libX = uibX if and only if X is uniform.

**Corollary 7.** Any uniform discrete set X is bounded up and down.

**Lemma 6.** If X and Y are discrete sets in $R$ and $X \subseteq Y$, then libY $\leq$ libX and uibY $\leq$ uibX.

**Proposition 6.** A function $f: X \to R$ defined on a bounded down discrete set X is $(q, r)$-continuous inside X for any $r \geq 0$ if $q <$ libX.

Proof. Indeed, in the discrete set $X$ with $q <$ lib$X$, there is $k > 0$ such that the only possibility that the inequality d($a$, $b$) $< q + k$ for two points $a$ and $b$ from $X$ is when $a = b$. This implies that conditions from Definition 3 are true for any $r \geq 0$.

**Corollary 8.** If $X = \{\ ku;\ k = 0, \pm 1, \pm 2, \pm 3,\ \ldots\ \}$ or $X = \{\ ku;\ k \in \mathbf{Z}$ and $m < k < n\ \}$ and $q < u$, then any function $f: \mathbf{R} \to \mathbf{R}$ is $(q, r)$-continuous inside $X$ for any $r > 0$.

**Corollary 9.** Any function $f: X \to \mathbf{R}$ defined only on a discrete set $X = \{\ ku;\ k = 0, \pm 1, \pm 2, \pm 3,\ \ldots\ \}$ or $X = \{\ ku;\ k \in \mathbf{Z}$ and $m < k < n\ \}$ is $(q, r)$-continuous if $q < u$.

**Proposition 7.** If $X$ and $Y$ are bounded up discrete sets, then $X \cup Y$ is also a bounded up discrete set.

Proof. For the union $X \cup Y$, there are two possibilities: either all points from one of these sets (say, $X$) are situated to the right of all points from the other of these sets (say, $Y$) or there are three points $x$, $z \in X$ and $y \in Y$ such that $y$ is situated between $x$ and $z$.

By the initial conditions, $X$ and $Y$ are bounded up discrete sets. It means that all distances between any two adjacent points in $X$ are less than some number $k > 0$ and all distances between any two adjacent points in $Y$ are less than some number $h > 0$. Then in the first of the two considered options, all distances between any two adjacent points in $X \cup Y$ are less than $q = \max\{k, h, r\}$ where $r$ is the distance between $X$ and $Y$. It means that $X \cup Y$ is also a bounded up discrete set.

In the case when there are three points $x$, $z \in X$ and $y \in Y$ such that $y$ is situated between $x$ and $z$, we have a different situation where for any two adjacent points $a$ and $b$, either they both belong to the same set $X$ or $Y$ or they belong to different sets but one of them lies between two points from the other set (say, $b$ lies between $a$ and $c$). In the first case, the distance between $a$ and $b$ is less than $p = \max\{k, h\}$. In the second case, the distance between $a$ and $b$ is less than or equal to (when $b = c$) the distance between $a$ and $c$. Consequently, the distance between $a$ and $b$ is less than $p = \max\{k, h\}$.

Proposition is proved.

### 3. Discrete and Fuzzy Intermediate Value Theorems

Let $f: [a, b] \to \mathbf{R}$ be an $r$-continuous function.

**Theorem 1 (the First Fuzzy Bolzano-Cauchy Theorem).** If $f(a) < 0$ and $f(b) > 0$, then there is [at least one] point $c \in [a, b]$ such that $|f(c)| < r$.

Proof. We consider several possibilities. First, there is some point $c \in [a, b]$, for which $f(c) = 0$. In this case, the statement of the theorem is true.

Second option is that there is no such a point, but the function $f(x)$ is larger (or less) than zero only on a finite set of points $x_0, x_1, \ldots, x_n$ from $[a, b]$. Then there is a sequence $l = \{ z_i; i = 1, 2, 3, \ldots \}$ of points from $[a, b]$ such that $\lim_{i \to \infty} z_i = c$ and $f(z_i) < 0$ for all $i = 1, 2, 3, \ldots$.

As $f(x)$ is an $r$-continuous function, it is bounded on $[a, b]$ (Burgin, 1995). As a closed interval is a compact space, it is possible to choose such a subsequence of the sequence $l = \{ z_i; i = 1, 2, 3, \ldots \}$ that the values of $f(x)$ on the elements of this subsequence have a limit. Thus, it is possible to assume that the sequence $f(l) = \{ f(z_i) ; i = 1, 2, 3, \ldots \}$ has a limit.

Let $\lim_{i \to \infty} f(z_i) = u$. Then by the definition of fuzzy continuity, we have $|u - f(x_0)| < r$. All numbers $f(z_i)$ are less than zero. Consequently, $u \leq 0$, while $f(c) > 0$. As a result, we have $u \leq 0 < f(c)$. This implies $|f(x_0)| = |0 - f(c)| \leq |u - f(c)| < r$. So, in this case, the theorem is also true.

The case when the function $f(x)$ is less than zero only on a finite set of points $x_0, x_1, \ldots, x_n$ from $[a, b]$ is considered in a similar way with the same result of validity the theorem.

The last option is when there are infinitely many points in which $f(x)$ is larger than zero and there are infinitely many points in which $f$ is less than zero. The we can take the subset $A = \{ x \in [a, b] ; f(x) > 0 \}$ and $B = \{ x \in [a, b] ; f(x) < 0 \}$. Then $[a, b] = A \cup B$. Any interval is a connected set. Consequently, there is such a point $x_0$ that belongs to the closures both of $A$ and of $B$. By the definition of a closure in $\mathbf{R}$, there are a sequence $l = \{ z_i; i = 1, 2, 3, \ldots \}$ of points from $B$ and a sequence $h = \{ x_i; i = 1, 2, 3, \ldots \}$ of points from $A$ such that $\lim_{i \to \infty} z_i = c$ and $\lim_{i \to \infty} x_i = c$.

As $f(x)$ is a bounded function on a closed interval $[a, b]$, it is possible to assume that both sequences $f(l) = \{ f(z_i) ; i = 1, 2, 3, \ldots \}$ and $f(h) = \{ f(x_i) ; i = 1, 2, 3, \ldots \}$ have limits.

Let $\lim_{i \to \infty} f(z_i) = u$ and $\lim_{i \to \infty} f(x_i) = v$. Then by the definition of fuzzy continuity and properties of the absolute value, we have $|u - f(c)| < r$ and $|v - f(c)| < r$. As we have $f(c) \neq 0$ (cf. case 1), then either $f(c) > 0$ or $f(c) < 0$.

Let us consider at first the first case. All numbers $f(z_i)$ are less than zero. Consequently, $u \leq 0$, while $f(c) > 0$. As a result, we have $u \leq 0 < f(c)$. This implies the following inequalities $|f(c)| = |0 - f(c)| \leq | u - f(c)| < r$. So in this case, the theorem is true.

In the second case, $f(c) < 0$ and all numbers $f(x_i)$ are larger than zero. Consequently, $v \geq 0$, while $f(c) < 0$. As a result, we have $f(c) < 0 \leq v$. This implies the following inequalities $|f(c)| = |0 - f(c)| \leq | v - f(c)| < r$. So in this case, the theorem is also true.

Theorem is proved because we have considered all possible cases.

Let $f: [a, b] \to R$ be a continuous function.

**Corollary 10 (the First Bolzano-Cauchy Theorem).** If $f(a) < 0$ and $f(b) > 0$, then there is [at least one] point $c \in [a, b]$ such that $|f(c)| = 0$.

By Corollary 4, any $(q, r)$-continuous function is $r$-continuous. Thus, Theorem 1 implies the following result.

**Corollary 11.** If $q \geq u$ and $v \geq r$, then for any $(q, r)$-continuous function $f: X \to Y$, and any $a$ and $b$ from $X$ with $f(a) < 0$ and $f(b) > 0$, there is an element $c$ from $X$ such that $|f(c)| < r$.

Let us consider two bounded up and down discrete sets $X = \{ a_k ; k = 0, \pm 1, \pm 2, \pm 3, \ldots \}$ or $X = \{ a_k ; k = 1, 2, 3, \ldots , n\}$ and $Y = \{ 0, b_k ; k = \pm 1, \pm 2, \pm 3, \ldots \}$. By Definition 7, there are positive real numbers $u$ and $v$ such that the distance between any two points in $X$ is less than or equal to $u = \text{uib}X$ and the distance between any two points in $Y$ is larger than or equal to $v = \text{lib}Y$.

**Theorem 2 (the First Discrete Bolzano-Cauchy Theorem).** If $q \geq u$ and $v \geq r$, then for any $(q, r)$-continuous function $f: X \to Y$, and any $a$ and $b$ from $X$ with $f(a) < 0$ and $f(b) > 0$, there is an element $c$ from $X$ such that $f(c) = 0$.

Proof. By Corollary 4, the function $f(x)$ is $r$-continuous. Let us extend this function to a function defined on the whole space $\mathbf{R}$ when $X$ is infinite and to a function defined in the interval $[c, d]$ when $X$ is finite and is a subset of this interval. Let us consider the function $g(x)$ equal to $f(a_k)$ when $a_k \leq x < a_{k+1}$ and $k = 0, \pm 1, \pm 2, \pm 3, \ldots$ in the first case or $k = 1, 2, 3, \ldots, n$ in the second case. This function provides such an extension. By Definition 3, $g(x)$ is $(q, r)$-continuous because $|g(x) - g(z)| \leq r$ whenever $|x - z| \leq q$ as $q \geq u = \text{uib}X$. Thus, by Corollary 4, the function $g(x)$ is $r$-continuous.

Then by Theorem 1, there is an element $c$ from the interval $[a, b]$ such that $a < c < b$ and $|g(c)| < r$. By the definition of the function $g(x)$, we have

$$g(c) = g(a_k) = f(a_k)$$

for some $i \in \mathbf{Z}$ and $a_i$ that belongs to the interval $[a, b]$.

Consequently, we have

$$|f(a_i)| < r \qquad (1)$$

As $f(a_i)$ belongs to the set $Y$, 0 belongs to the set $Y$, and the distance between points in $Y$ is larger than or equal to $r$ (because $v = \text{lib}Y \geq r$), the inequality (1) is possible only if $f(a_i) = 0$. Thus, we can take $c = a_i$ and have $a < c < b$ and $f(c) = 0$.

Theorem is proved.

Let us consider two discrete sets $X = \{ku; k = 0, \pm 1, \pm 2, \pm 3, \ldots\}$ (or $X = \{ku; k \in \mathbf{Z} \text{ and } m < k < n\}$) and $Y = \{kv; k = 0, \pm 1, \pm 2, \pm 3, \ldots\}$ where $u$ and $v$ are positive real numbers. They are uniform and thus, bounded down discrete sets. Then Theorem 2 implies the following result.

**Corollary 12.** If $q \geq u$ and $v \geq r$, then for any $(q, r)$-continuous function $f: X \to Y$, and any $a = k_1 u$, $b = k_2 v$ with $f(a) < 0$ and $f(b) > 0$, there is an element $c$ from $X$ such that $f(c) = 0$.

**Corollary 13 (the First Uniform Discrete Bolzano-Cauchy Theorem).** If $X$ and $Y$ are uniform discrete spaces (grids) with grid spacings $u$ and $v$, correspondingly, $q \geq$

$u$ and $v \geq r$, then for any $(q, r)$-continuous function $f: X \to Y$, and any $a$ and $b$ from $X$ with $f(a) < 0$ and $f(b) > 0$, there is an element $c$ from $X$ such that $f(c) = 0$.

Let $f: [a, b] \to \mathbf{R}$ be an $r$-continuous function.

**Theorem 3 (the Fuzzy Intermediate Value Theorem).** If $f(a) = k$, $f(b) = h$, and $k < h$ ($k > h$), then for any number $l$ from the interval $[k, h]$ (from the interval $[h, k]$) there is [at least one] point $c \in [a, b]$ such that $|l - f(c)| < r$.

Proof. We prove the theorem for the case when $k < h$ because the case when $k > h$ is proved in a similar way. Let us take a number $l$ from the interval $[k, h]$. If $l = k$, then $f(a) = l$ and we have the necessary result. If $l = h$, then $f(b) = l$ and we have the necessary result.

If $k < l < h$, we consider the function $g(x) = f(x) - l$. For this function, we have $g(a) < 0$ and $g(b) > 0$. Besides, as the difference of an $r$-continuous function and continuous function is an $r$-continuous function (Burgin, 1995), $g: [a, b] \to \mathbf{R}$ is also an $r$-continuous function. By Theorem 1, there is such a point $c \in [a, b]$ that $|0 - g(c)| < r$. The natural metric on $\mathbf{R}$ is invariant with respect to translations, i.e., $|a - b| = |(a + c) - (b + c)|$. Consequently, $|l - f(c)| = |l - g(c) + l| = |0 - g(c)| < r$.

Theorem is proved.

For the case when $r = 0$, we have the following classical result.

Such a classical result as the Intermediate Value Theorem is a direct corollary of Theorem 3.

By Corollary 4, any $(q, r)$-continuous function is $r$-continuous. Thus, Theorem 3 implies the following result.

**Corollary 14.** If a function $f: X \to Y$ is $(q, r)$-continuous inside the interval $[a, b]$, $f(a) = k$, $f(b) = h$, and $k < h$ ($k > h$), then for any number $l$ from the interval $[k, h]$ (from the interval $[h, k]$) there is [at least one] point $c \in [a, b]$ such that $|l - f(c)| < r$.

Let us look what properties of functions on discrete sets can be deduced from Theorem 3.

Such a classical result as the Intermediate Value Theorem is also a corollary of Theorem 3.

Let $f: [a, b] \to \mathbf{R}$ be a continuous function. Then we have the following classical result (cf., for example, (Ross, 1996)).

**Corollary 15 (the Intermediate Value Theorem).** If $f(a) = k$, $f(b) = h$, and $k < h$ ( $k > h$ ), then for any number $l$ from the interval $[k, h]$ ( from the interval $[h, k]$ ) there is [at least one] point $c \in [a, b]$ such that $f(c) = l$.

Let us consider two bounded up and down discrete sets $X = \{\ a_k\ ;\ k = 0, \pm 1, \pm 2, \pm 3, \ldots\ \}$ or $X = \{\ a_k\ ;\ k = 1, 2, 3, \ldots, n\}$ and $Y = \{\ 0, b_k\ ;\ k = \pm 1, \pm 2, \pm 3, \ldots\ \}$. By Definition 7, there are positive real numbers $u$ and $v$ such that the distance between any two points in $X$ is less than or equal to $u = \text{uib}X$ and the distance between any two points in $Y$ is larger than or equal to $v = \text{lib}Y$.

**Theorem 4 (the Discrete Intermediate Value Theorem).** If $q \geq u$ and $v \geq r$, then for any $(q, r)$-continuous function $f: X \to Y$, any $a$ and $b$ from $X$ and any element $l$ that belongs to the intersection $Y \cap [f(a), f(b)]$ (or to the intersection $Y \cap [f(b), f(a)]$ when $f(a) > f(b)$ ), there is [at least one] point $c$ from $X$ such that $a < c < b$ and $f(c) = l$.

<u>Proof</u>. By Corollary 4, the function $f(x)$ is $r$-continuous. Let us extend this function to a function defined on the whole space $\mathbf{R}$ when $X$ is infinite and to a function defined in the interval $[c, d]$ when $X$ is finite and is a subset of this interval. Let us consider the function $g(x)$ equal to $f(a_k)$ when $a_k \leq x < a_{k+1}$ and $k = 0, \pm 1, \pm 2, \pm 3, \ldots$ in the first case or $k = 1, 2, 3, \ldots, n$ in the second case. This function provides such an extension. By Definition 3, $g(x)$ is $(q, r)$-continuous because $|\,g(x) - g(z)\,| \leq r$ whenever $|\,x - z\,| \leq q$ as $q \geq u = \text{uib}X$. Thus, by Corollary 4, the function $g(x)$ is $r$-continuous.

Then by Theorem 1, there is an element $c$ from the interval $[a, b]$ such that $a < c < b$ and $|l - g(c)| < r$. By the definition of the function $g(x)$, we have

$$g(c) = g(a_k) = f(a_k)$$

for some $i \in \mathbf{Z}$ and $a_i$ that belongs to the interval $[a, b]$.

Thus, we have

$$|l - f(ku)| < r \qquad (2)$$

As points $l$ and $f(a_i)$ both belong to the set $Y$, and the distance between points in $Y$ is larger than or equal to $r$ (because $v = \text{lib} Y \geq r$), the inequality (2) is possible only if $f(a_i) = l$. Thus, we can take $c = a_i$ and have $a < c < b$ and $f(c) = l$.

Theorem is proved.

Let us consider two discrete sets $X = \{\ ku;\ k = 0, \pm 1, \pm 2, \pm 3, \ldots\ \}$ (or $X = \{\ ku;\ k \in \mathbf{Z}$ and $m < k < n\ \}$) and $Y = \{\ kv;\ k = 0, \pm 1, \pm 2, \pm 3, \ldots\ \}$.

**Corollary 16.** If $q \geq u$ and $v \geq r$, then for any $(q, r)$-continuous function $f\colon X \to Y$, any $a = k_1 u$, $b = k_2 v$ with $k_1 < k_2$ ($m < k_1 < k_2 < n$) and any element $l$ that belongs to the intersection $Y \cap [f(a), f(b)]$ (or to the intersection $Y \cap [f(b), f(a)]$ when $f(a) > f(b)$ ), there is [at least one] point $c$ from $X$ such that $a < c < b$ and $l = f(c)$.

**Corollary 17.** If a function $f\colon X \to Y$ is $(u, v)$-continuous inside the interval $[a, b]$ with $a = k_1 u$, $b = k_2 v$ and $k_1 < k_2$ ($m < k_1 < k_2 < n$), then for any element $l$ that belongs both to the set $Y$ and to the interval $[f(a), f(b)]$ (or to the interval $[f(a), f(b)]$ when $f(a) > f(b)$ ), there is an element $c$ from $X$ such that $a < c < b$ and $l = f(c)$.

**Corollary 18 (the Discrete Uniform Intermediate Value Theorem).** If $X$ and $Y$ are uniform discrete spaces (grids) with grid spacings $u$ and $v$, correspondingly, $q \geq u$ and $v \geq r$, then for any $(q, r)$-continuous function $f\colon X \to Y$, any $a$ and $b$ from $X$ and any element $l$ that belongs to the intersection $Y \cap [f(a), f(b)]$ (or to the intersection $Y \cap [f(b), f(a)]$ when $f(a) > f(b)$ ), there is [at least one] point $c$ from $X$ such that $a < c < b$ and $f(c) = l$.

This implies the Intermediate Value Theorem that was obtained in digital topology for the case $u = v = 1$ (cf. (Rosenfeld, 1979; 1986; Hamlet, 2002)). Namely, taking a discrete set $X$ that is either the set of all integer numbers $\mathbf{Z}$ or an interval in $\mathbf{Z}$, we have the following result.

**Theorem 5 (Digital Topology Intermediate Value Theorem).** For any $(1, 1)$-continuous function $f\colon X \to \mathbf{Z}$, any $m$ and $n$ from $X$ and any element $l$ such that $m < n$ and $f(m) < l < f(n)$ or $f(m) > l > f(n)$, there is an element $c$ from $X$ such that $m < c < n$ and $l = f(c)$.

## 4. Applications

The Intermediate Value Theorem is used to study continuous images of connected sets. In a similar way, the Fuzzy Intermediate Value Theorem can be used to study fuzzy continuous images of fuzzy connected sets introduced and studied in (Burgin, 2004). To simplify exposition, here we consider fuzzy connected sets not in general topological spaces as in (Burgin, 2004), but only in the real line.

Distance between two points in a real line is defined as $d(a, c) = |a - c|$. It is possible to extend this distance to distance between a set and a point.

**Definition 9** (Kuratowski, 1968). If $A \subseteq R$ and $c \in R$, then the distance between $A$ and $c$ is defined as

$$\text{dist}(c, A) = \text{dist}(A, c) = \inf \{ d(a, c) ; a \in A \}.$$

**Remark 2.** If $d(A, c) > 0$, then $c \notin A$.

**Definition 10.** a) A subset $C$ of $R$ is called *r-disconnected* if $C = A \cup B$, for any point $a \in A$, the distance $\text{dist}(a, B) > r$ and for any point $b \in B$, the distance $\text{dist}(A, b) > r$. Sets $A$ and $B$ are called *mutually r-disconnected*.

b) A subset $C$ of $X$ is called *r-connected* in $X$ if it is not *r-disconnected*.

When the structure Q is fixed, all Q-connected (Q-disconnected) in $X$ sets are called fuzzy or relatively connected (disconnected) in $X$.

**Example 7.** Points 0 and 0.7 form a 1-connected set, while points 0 and 1.1 or points 1 and 4 form a 1-disconnected set.

**Example 8.** Let us consider sets $A = [0,1] \cup [1.7, 3]$, and $B = [0,1] \cup [2.6, 3]$. The set $A$ is 1-connected, while the set $B$ is 1-disconnected.

Directly from Definition 9, we obtain the following results.

**Lemma 7.** A set $C$ is connected (disconnected) if and only if it is 0-connected (0-disconnected).

**Lemma 8.** a) If a set $C$ is *r-connected*, then it is *q-connected* for any $q > r$.

b) If a set $C$ is *r-disconnected*, then it is *q-disconnected* for any $q < r$.

**Corollary 17.** Any connected set is *r-connected* for any $r > 0$.

**Lemma 9.** If sets $A$ and $B$ are called mutually $r$-disconnected for some $r \in \mathbf{R}^+$, then $A \cap B = \emptyset$.

These results demonstrate that the concept of $r$-connectedness is a *natural extension* of the concept of conventional connectedness.

Let $q, r \in \mathbf{R}^+$.

**Theorem 6.** If $C$ is a $q$-connected set and a real function $f(x)$ is $(q, r)$-continuous, then $f(C)$ is an $r$-connected set.

Proof. Let us assume that $C$ is a $q$-connected set, a function $f(x)$ is $(q, r)$-continuous but the set $D = f(C)$ is a $r$-disconnected. By Definition 10, $D = A \cup B$ and for any point $a \in A$, the distance $\text{dist}(a, B) > r$ and for any point $b \in B$, the distance $\text{dist}(A, b) > r$.

Let us consider sets $E = f^{-1}(A)$ and $H = f^{-1}(B)$. By Lemma 9, $A \cap B = \emptyset$. As $f(x)$ is a function, $E \cap H = \emptyset$. As $C$ is a $q$-connected set, there is a point $c \in E$ such that the distance $\text{dist}(c, H)$ is smaller than or equal to $q$. Then $f(c) = a \in A$.

The distance $\text{dist}(a, B) > r$ because $A$ and $B$ are $r$-disconnected sets. Consequently, there is $k > 0$ such that $\text{dist}(a, B) > r + k$. By Definition 9, $d(a, b) > r + k$ for any point $b \in B$. At the same time, the function $f(x)$ is $(q, r)$-continuous. Thus, by Definition 3, for any $\varepsilon > 0$ there is $\delta > 0$ such that the inequality $| a - x | < q + \delta$ implies the inequality $| f(x) - f(a) | < r + \varepsilon$. Consequently, there is $\delta > 0$ such that the inequality $| a - x | < q + \delta$ implies the inequality $| f(x) - f(a) | < r + k$.

By the choice of the point $c$, we have $\text{dist}(c, H) \leq q$. Thus by Definition 9, there is a point $d \in H$ such that the distance $d(c, d) < q + \delta$. Consequently, we have $| f(c) - f(d) | < r + k$ and $f(d) \in B$. However, as we know, $d(a, b) > r + k$ for any point $b \in B$. This contradiction shows that the set $D = f(C)$ is not $r$-disconnected. Thus, the set $D$ is a $r$-connected.

Theorem is proved.

**Corollary 18.** If $D$ is an $r$-disconnected set and a real function $f(x)$ is $(q, r)$-continuous, then $f^{-1}(D)$ is a $q$-disconnected set.

Theorem 6 also implies the following classical result (cf., for example, (Dieudonné, 1960) or (Kuratowski, 1968)).

**Corollary 19.** If $C$ is a connected set and a real function $f(x)$ is continuous, then $f(C)$ is a connected set.

The Intermediate Value Theorem is often used to prove existence of solutions for different kinds of equations. For instance, existence of roots of an odd degree polynomial is a direct corollary of this theorem. In a similar way, it is possible to use the Discrete Intermediate Value Theorem for proving existence of computable with given precision solutions for different kinds of equations.

The Intermediate Value Theorem is often used to find inverse functions and their properties. Here we use the Fuzzy Intermediate Value Theorem and Discrete Intermediate Value Theorem for the same purposes.

**Definition 11.** A *gap* in a set $X$ of real numbers is an interval $I$ of real numbers that does not contain point from $X$ and there is no other interval that contains $I$ and has the same property.

**Proposition 4.** A set $X$ of real numbers is discrete if and only if any its point is adjacent to a gap from the right and to a gap from the left.

Let us take a strictly increasing (decreasing) real function $f(x)$ defined on an interval $I$ in $\mathbf{R}$. The interval $I$ may be closed or not, finite or infinite.

**Theorem 7.** If the length of any gap in the image $f(I)$ of the interval $I$ is less than $r$, then $f(x)$ is an $r$-continuous function.

<u>Proof.</u> We consider, at first, a strictly increasing real function $f(x)$ and show that $f(x)$ satisfies the condition from Definition 1.

Let us consider a point $a$ from the interval $I$ such that $a$ is not the right endpoint of $I$ and some $\varepsilon > 0$. We can take $\varepsilon$ so small that $y_1 = y + \varepsilon$ is less than the right endpoint of $f(I)$ if this endpoint exists where $y = f(a)$. As any gap in the image $f(I)$ of the interval $I$ is less than $r$, there is a point $y_2$ in the image $f(I)$ such that $y_1 \leq y_2 < y_1 + r$ and $y_2 = f(x_2)$ for some $x_2$ from the interval $I$.

As $y < y_1 \leq y_2$, we have $a < x_2$. Let us take $\delta = x_2 - a$ and consider a point $x$ from the interval $I$ such that $a < x < x_2 = a + \delta$. As $f(x)$ is a strictly increasing real function, we have $f(a) < f(x) < f(x_2) < y_1 + r = y + r + \varepsilon = f(a) + r + \varepsilon$. As $x$ is an arbitrary point between *a* and *a* + δ, this means that the function *f*(*x*) is right *r*-continuous at the point *a*.

In a similar way, we can show that the function *f*(*x*) is left *r*-continuous at the point *a*. Thus, the function *f*(*x*) is *r*-continuous at the point *a*.

The case of a strictly decreasing real function *f*(*x*) is treated in a similar way.

Theorem is proved.

Theorem 7 implies the following classical result.

**Corollary 20.** If the image *f*(*I*) of the interval *I* is also an interval, then *f*(*x*) is a continuous function.

By a similar technique, we can prove the following result.

Let us take a strictly increasing (decreasing) real function *f*(*x*) defined on a subset *X* of an interval *I* in **R**.

**Theorem 8.** If the length of any gap in *X* is less than *q* and the length of any gap in the image *f*(*X*) of *X* is less than *r*, then *f*(*x*) is an *r*-continuous function.

## 5. Conclusion

Thus, we have proved different extensions of the Intermediate Value Theorem: the Fuzzy Intermediate Value Theorem, Discrete Intermediate Value Theorem and Digital Topology Intermediate Value Theorem. These theorems have been applied to study images of connected and fuzzy connected sets, as well as to prove existence and fuzzy continuity of inverse functions.

The Intermediate Value Theorem is often used to prove existence of solutions for different kinds of equations. For instance, existence of roots of an odd degree polynomial is a direct corollary of this theorem. In a similar way, it is possible to use the Discrete Intermediate Value Theorem for proving existence of computable with given precision solutions for different kinds of equations.

It would be interesting to study similar properties of fuzzy continuous functions in multidimensional Euclidean spaces, as well as in general metric spaces in a similar way as, for example, Dieudonné (1960) studies continuous functions in general metric spaces.

It is necessary to remark that this research also opens new venues for scalable topology (Burgin, 2004; 2005; 2006). For instance, it is possible to study discrete sets in scalable topological spaces.

**R e f e r e n c e s**

1. Boxer, L. (1994) Digitally continuous functions, *Pattern Recognition Letters*, v. 15, pp. 833-839
2. Burgin, M. (1995) Neoclassical Analysis: Fuzzy Continuity and Convergence, *Fuzzy Sets and Systems*, v. 75,  pp. 291-299
3. Burgin, M. (1999) General Approach to Continuity Measures,  *Fuzzy Sets and Systems*, v. 105, No. 2, pp. 225-231
4. Burgin, M. *Monotonicity, Fuzzy Extrema, and Fuzzy Conditional Derivatives of Real Functions*, University of California, Los Angeles, Mathematics Report Series, MRS Report 03-13, 2003, 73 p.
5. Burgin, M. (2004) *Discontinuity Structures in Topological Spaces*, International Journal of Pure and Applied Mathematics, 2004, v. 16, No. 4, pp. 485-513
6. Burgin, M. (2006) Scalable Topological Spaces, 5$^{th}$ *Annual International Conference on Statistics, Mathematics and Related Fields*, 2006 Conference Proceedings, Honolulu, Hawaii, January, pp. 1865-1896
7. Burgin, M. *Fuzzy Continuity in Scalable Topology*, Preprint in Mathematics, math/0512627 (subjects: math.GN; math-ph), 2005, 30 p.    (electronic edition: http://arXiv.org)
8. Burgin, M. *Neoclassical Analysis*: *Calculus closer to the Real World*, Nova Science Publishers, Hauppauge, N.Y., 2008
9. Cartwrite, J.H.E., and Piro, O. (1992) The Dynamics of Runge-Kutt Methods, Int. Journ. of Bifurcation and Chaos, v. 2, No. 3, pp. 427-450
10. Dieudonné, J. *Foundations of Modern Analysis*, Academic Press, New York and London, 1960
11. Gontar, V. (1997) Theoretical Foundation for the Discrete Dynamics of Physicochemical Systems: Chaos, Self-Organization, Time and Space in Complex Systems, Discrete Dynamics in Nature and Society, v. 1, No. 1, pp.31-43
12. Hamlet, D. (2002) Continuity in Software Systems, in *Proceedings of the ISSTA* 02, Rome, pp. 196-200
13. Hänel, D., Roth, P., Rose, M.,  Thill, C., Uphoff U. and Vilsmeier R. (1999) Adaptive Grid Methods for Reactive Flows, in *Hyperbolic problems*: *theory, numerics, applications*, Birkhäuser, pp. 445-454
14. Kuratowski, K. *Topology*, Academic Press, Waszawa, v. 1, 1966; v. 2, 1968

15. Rosenfeld, A. (1979) Digital topology, *American Mathematical Monthly*, v. 86, pp. 621-630.

16. Rosenfeld, A. (1986) Continuous functions on digital pictures, *Pattern Recognition Letters*, v. 4, pp. 177-184.

17. Ross, K.A. *Elementary Analysis*: *The Theory of Calculus*, Springer-Verlag, New York/Berlin/Heidelberg, 1996